\documentclass[a4paper,11pt]{article}
\usepackage[french]{babel}
\usepackage{amsmath,amscd,amssymb}

\input xy
\xyoption{all}
\usepackage[all]{xy}

\newcommand{\PP}{{\bf{P}}}
\newcommand{\CC}{{\bf{C}}}
\newcommand{\ZZ}{{\bf{Z}}}

\newcommand{\GG}{{\bf{G}}}

\newcommand{\corps}{\textup{k}}
\newcommand{\V}{\textup{V}}
\newcommand{\B}{\textup{B}}
\newcommand{\X}{\textup{X}}
\newcommand{\Y}{\textup{Y}}
\newcommand{\Z}{\textup{Z}}
\newcommand{\K}{\textup{K}}

\newcommand{\C}{\textup{C}}
\newcommand{\N}{\textup{N}}
\newcommand{\F}{\textup{F}}
\newcommand{\U}{\textup{U}}
\newcommand{\E}{\textup{E}}
\newcommand{\G}{\textup{G}}
\newcommand{\I}{\textup{I}}
\newcommand{\J}{\textup{J}}
\newcommand{\A}{\textup{A}}
\newcommand{\D}{\textup{D}}
\newcommand{\T}{\textup{T}}
\newcommand{\M}{\textup{M}}
\newcommand{\coh}{\textup{H}}
\newcommand{\lb}{\textup{L}}
\newcommand{\sur}{\textup{S}}

\newcommand{\der}[1]{\Theta_{#1}}

\newcommand{\re}[2]{(\textit{voir} \cite{#1}, #2)}
\newcommand{\res}[1]{(\textit{voir} \cite{#1})}

\newcommand{\qed}{\hfill$\blacksquare$}

\newcounter{compteur}
\newcounter{sect}
\newcommand{\sect}[1]{\\$\ $\newline\stepcounter{sect}\textbf{\thesect}.
\textbf{#1}\setcounter{compteur}{1}\\$\ $\newline}

\newcommand{\lemme}[1]{\\$\ $\newline\textbf{Lemme
\thesect.\thecompteur}.$-$\noindent\textit{#1}\stepcounter{compteur}\\
$\ $\newline}

\newcommand{\lemmeref}[3]{\\$\ $\newline\textbf{Lemme
\thesect.\thecompteur} (\textit{voir} \cite{#1}, #2).$-$\noindent
\textit{#3}\stepcounter{compteur}\\$\ $\newline}

\newcommand{\prop}[1]{\\$\ $\newline\textbf{Proposition
\thesect.\thecompteur}.$-$\noindent\textit{#1}\stepcounter{compteur}\\
$\ $\newline}

\newcommand{\marque}{(\textbf{\thesect.\thecompteur) }\stepcounter{compteur}}

\newcommand{\rem}{\\$\ $\newline\textbf{Remarque \thesect.\thecompteur}.$-$\noindent
\stepcounter{compteur}}

\newcommand{\ex}{\\$\ $\newline\textbf{Exemple \thesect.\thecompteur}.$-$\noindent
\stepcounter{compteur}}

\begin{document}

\centerline{\large{\textbf{Caract\'erisation de l'espace projectif}}}
$\ $
\vspace{0.5cm}\\
$\ $
\centerline{St\'ephane \textsc{Druel}}
$\ $
\vspace{1cm}\\
\noindent\textbf{Introduction}\\
\newline
\indent Soit $\X$ une vari\'et\'e projective \textit{normale}, de dimension 
$n\ge 2$, d\'efinie 
sur un corps 
alg\'ebriquement clos de caract\'eristique nulle.
Un \textit{feuilletage en courbes sur $\X$}
est la donn\'ee d'un fibr\'e en droites $\lb$ sur $\X$ et
d'une application non nulle $\eta : \Omega_{\X}^{1}\longrightarrow\lb^{-1}$.\\
\indent Le point de d\'epart de notre travail 
est la caract\'erisation suivante de l'espace projectif $\PP^{n}$. 
Si
$\lb$ est \textit{ample} alors $\lb\simeq\mathcal{O}_{\X}(\Y)$, o\`u $\Y\subset\X$ est un 
diviseur irr\'eductible et normal, et $\X$ s'identifie \`a la contraction de la section de 
$\PP_{\Y}(\mathcal{O}_{\Y}\oplus\mathcal{O}_{\Y}(-\Y))\longrightarrow\Y$ de fibr\'e
normal $\mathcal{O}_{\Y}(-\Y)$ \re{Wa83}{Theorem 2}. 
Si, de plus, $\X$ est lisse alors $\X$ est isomorphe \`a l'espace projectif 
$\PP^{n}$ \re{Wa83}{Theorem 1}.\\
\indent Les techniques utilis\'ees par Wahl rel\`event de l'alg\`ebre
commutative. Nous donnons une d\'emonstration nouvelle, de nature g\'eom\'etrique,
de son r\'esultat et d'un \'enonc\'e plus g\'en\'eral.\\
\indent Le $\textit{lieu singulier}$ du feuilletage $(\eta,\lb)$ est le 
sous-sch\'ema ferm\'e $\Z(\eta,\lb)$ de $\X$ dont l'id\'eal est l'image de 
l'application induite 
$\Omega_{\X}^{1}\otimes\lb\longrightarrow\mathcal{O}_{\X}$.
Une courbe irr\'eductible $\C\subset\X$ est appel\'ee une $\textit{feuille}$ si 
$\C\not\subset\Z(\eta,\lb)$ et si la restriction $\eta_{|C}$ de $\eta$ 
\`a $\C$ se factorise \`a travers l'application naturelle
${\Omega_{\X}^{1}}_{|C}\longrightarrow\Omega_{\C}^{1}$.\\
$\ $
\newline
\noindent\textbf{Th\'eor\`eme}.$-$\textit{Soient 
$\X$ une vari\'et\'e projective normale, de dimension $\ge 2$, d\'efinie 
sur un corps alg\'ebriquement clos de caract\'eristique nulle, 
et $(\eta,\lb)$ un feuilletage en courbes sur $\X$. Si $\lb\cdot\C>0$ pour toute courbe $\C\subset\X$
et si $\Z(\eta,\lb)\neq\emptyset$ alors $\lb\simeq\mathcal{O}_{\X}(\Y)$, o\`u $\Y\subset\X$ est un 
diviseur irr\'eductible et normal, et
$\X$ s'identifie \`a la contraction de la section du morphisme
$\PP_{\Y}(\mathcal{O}_{\Y}\oplus\mathcal{O}_{\Y}(-\Y))\longrightarrow\Y$ de fibr\'e
normal $\mathcal{O}_{\Y}(-\Y)$.}\\
$\ $
\newline
\indent La preuve de ce r\'esultat repose sur un crit\`ere d'alg\'ebricit\'e
des feuilles d'un feuilletage alg\'ebrique d\'emontr\'e de mani\`ere ind\'ependante
par Bost \re{Bo01}{Theorem 3.5} et Bogomolov-McQuillan \res{BoMc01}. Si 
$\Z(\eta,\lb)\neq\emptyset$, nous montrons que les feuilles du feuilletage $(\eta,\lb)$ 
sont des courbes rationnelles passant par un m\^eme point de $\X$.\\
$\ $
\newline
\noindent\textbf{Corollaire}.$-$\textit{Soient 
$\X$ une vari\'et\'e projective lisse de dimension $n\ge 2$, d\'efinie 
sur un corps alg\'ebriquement clos de caract\'eristique nulle,
et $(\eta,\lb)$ un feuilletage en courbes sur $\X$. Si 
$\lb\cdot\C>0$ pour toute courbe $\C\subset\X$, alors 
ou bien $\X$ poss\`ede une fibration en droites projectives, ou bien
$\X$ est isomorphe \`a l'espace projectif $\PP^{n}$.}
\\
$\ $
\newline
\indent Le dernier r\'esultat est le suivant.\\
$\ $
\newline
\noindent\textbf{Proposition}.$-$\textit{Soient 
$\X$ une surface projective normale, d\'efinie 
sur un corps alg\'ebriquement clos de caract\'eristique nulle, 
et $(\eta,\lb)$ un feuilletage en courbes sur $\X$. 
Si $\lb\cdot\C>0$ pour toute courbe $\C\subset\X$, 
alors ou bien $\X$ est une surface g\'eom\'etriquement
r\'egl\'ee, ou bien $\X$ v\'erifie les conclusions du th\'eor\`eme.}
\sect{Feuilletage}
\marque Soit $\X$ une vari\'et\'e alg\'ebrique
d\'efinie sur un corps $\corps$. Un \textit{feuilletage en courbes sur $\X$}
est la donn\'ee d'un fibr\'e en droites $\lb$ sur $\X$ et
d'une application non nulle $\eta : \Omega_{\X}^{1}\longrightarrow\lb^{-1}$.
Son $\textit{lieu singulier}$ est le 
sous-sch\'ema ferm\'e $\Z(\eta,\lb)$ de $\X$ dont l'id\'eal est l'image de 
l'application induite 
$\Omega_{\X}^{1}\otimes\lb\longrightarrow\mathcal{O}_{\X}$.
Enfin, une courbe irr\'eductible $\C\subset\X$ est appell\'ee une $\textit{feuille}$ si 
$\C\not\subset\Z(\eta,\lb)$ et si la restriction $\eta_{|C}$ de $\eta$ 
\`a $\C$ se factorise \`a travers l'application naturelle
${\Omega_{\X}^{1}}_{|C}\longrightarrow\Omega_{\C}^{1}$, \textit{i.e.}, s'il existe
un diagramme commutatif\\
\centerline{
\xymatrix{ 
{{\Omega_{\X}^{1}}_{|C}} \ar[r]\ar[d]_-{\eta_{|C}} & {\Omega_{\C}^{1}}
\ar@/^/[dl]\\  
{{\lb^{-1}}_{|\C}}
}
}
\lemme{Soient $\X$ une vari\'et\'e alg\'ebrique
d\'efinie sur un corps $\corps$ de caract\'eristique nulle et 
$(\eta,\lb)$ un feuilletage en courbes sur $\X$. 
Soit enfin $n:\X^{n}\longrightarrow X$ la normalisation de $\X$. 
Il existe alors une unique application 
$\eta_{n}:\Omega_{\X^{n}}^{1}\longrightarrow n^{*}\lb^{-1}$
factorisant $n^{*}\eta : n^{*}\Omega_{\X}^{1}\longrightarrow n^{*}\lb^{-1}$ et, de plus,
$n^{-1}({\Z(\eta,\lb)}_{\textup{red}})\subset
{\Z(\eta_{n},n^{*}\lb)}_{\textup{red}}$.}
\textit{D\'emonstration}.$-$Soit $\textup{Spec}(\A)\subset\X$ un ouvert affine
au dessus duquel le fibr\'e $\lb$ est trivial. 
Soit $\textup{\K}$ le corps des fractions de $\A$ et
$\bar{\A}$ la cl\^oture int\'egrale de $\A$ dans $\textup{\K}$. Soit 
$\D_{\A}$ la d\'erivation de $\A$ d\'efinissant le feuilletage sur l'ouvert consid\'er\'e
et soit $\D_{\K}$ la d\'erivation de $\K$ correspondante. La d\'erivation
$\D_{\K}$ satisfait $\D_{\K}(\bar{\A})\subset\bar{\A}$ \res{Se66}
et c'est l'unique d\'erivation de $\bar{\A}$ induisant
$\D_{A}$ sur $\A$. La premi\`ere assertion du lemme est alors imm\'ediate.\\
\indent Soit $\I$ l'id\'eal de $\Z(\eta,\lb)\cap\textup{Spec}(\A)$ dans $\A$,
et soit $\mathfrak{m}\supset\I$ un id\'eal maximal. Soit $\J:=\mathfrak{m}\bar{A}$. 
Comme $\D_{\A}(\A)\subset\I$, $\D_{\K}(\J)\subset\J$
et $\D_{\K}(\sqrt{J})\subset\sqrt{J}$. La d\'erivation $\D_{\K}$ induit donc
une d\'erivation de l'anneau $\bar{\A}/\sqrt{J}$. Le sch\'ema 
$\textup{Spec}(\bar{\A}/\sqrt{J})$ est lisse sur $\corps$, la d\'erivation
pr\'ec\'edente est donc identiquement nulle et $\D_{\K}(\bar{\A})\subset\sqrt{J}$.\qed
\ex Si le corps de base 
est $\CC$ et si $\B$ est une courbe projective, lisse, 
de genre au moins 2, il existe un fibr\'e vectoriel 
$\G$ de rang 2 et de degr\'e $0$ sur $\B$ dont toutes les puissances sym\'etriques 
sont \textit{stables}. Soit $\sur=\PP_{\B}(\G)$ la surface r\'egl\'ee associ\'ee. 
Le fibr\'e tangent relatif $\T_{\sur/\B}\subset\T_{\sur}$ est alors
de degr\'e $>0$ sur toutes les courbes de $\sur$ \re{Ha70}{Example 10.6}.
\sect{Feuilletage sur les surfaces}
\indent Le r\'esultat de cette section est la proposition suivante.
\prop{Soient 
$\sur$ une surface normale, $\B$ une courbe projective lisse et 
$p_{\sur} : \sur\longrightarrow\B$ un morphisme \`a fibres connexes. 
Soit $\F$ une fibre g\'en\'erale de $p_{\sur}$.
Soit $\sur'$ une
surface alg\'ebrique et $q_{\sur} : \sur\longrightarrow\sur'$ un morphisme dominant. 
Soit $(\eta_{\sur},\lb_{\sur})$ un feuilletage en courbes sur $\sur$. 
Supposons le fibr\'e $\lb_{\sur}$ num\'eriquement effectif et supposons
qu'une courbe $\C\subset\sur$ est contract\'ee par $q_{\sur}$ si et seulement si
$\lb_{\sur}\cdot\C=0$. Supposons enfin que les courbes trac\'ees sur $\sur$ et
contract\'ees par $q_{\sur}$ sont horizontales relativement \`a $p_{\sur}$.
Sous ces hypoth\`eses, la surface $\sur\longrightarrow\B$ est g\'eom\'etriquement 
r\'egl\'ee, le feuilletage vertical et $\lb_{\sur}\cdot\F\in\{1,2\}$.
Si $\lb_{\sur}\cdot\F=2$ alors le feuilletage est r\'egulier et le morphisme $q_{\sur}$
fini.
Si $\lb_{\sur}\cdot\F=1$, alors le lieu singulier du feuilletage est l'unique 
section de $p_{\sur}$ contract\'ee par $q_{\sur}$.}
\textit{D\'emonstration}.$-$L'application compos\'ee
$$\N^{*}_{\F/\sur}\simeq\mathcal{O}_{\F}\longrightarrow{\Omega_{\sur}^{1}}_{|\F}
\longrightarrow (\lb_{\sur}^{-1})_{|\F}$$
est identiquement nulle car $\lb_{\sur}\cdot\F>0$ et il existe donc
une factorisation $\Omega_{\F}^{1}\longrightarrow (\lb_{\sur}^{-1})_{|\F}$.
Les fibres g\'en\'erales de $p_{\sur}$ sont donc des courbes rationnelles lisses et
$\lb_{\sur}\cdot\F\in\{1,2\}$.\\
$\ $
\newline
\noindent $1^{\textup{er}}$ cas : $\lb_{\sur}\cdot\F=1$.$-$Les fibres de $p_{\sur}$ sont
irr\'eductibles et g\'en\'eriquement r\'eduites. La surface $\sur\longrightarrow\B$
est donc g\'eom\'etriquement r\'egl\'ee.\\
\indent Soit $\G$ un fibr\'e de rang 2
sur $\B$, \textit{normalis\'e}, de sorte
que $\sur=\PP_{\B}(\G)$ \res{Ha77}. Notons $\C_{0}\subset\sur$ une section de $p_{\sur}$
telle que $\mathcal{O}_{\sur}(\C_{0})$ soit isomorphe au fibr\'e tautologique
$\mathcal{O}_{\sur}(1)$. Le groupe
$\textup{Num}(\sur)$ est libre de rang $2$ et engendr\'e par les classes de 
$\C_{0}$ et $\F$. Le produit d'intersection est donn\'e par les formules
$\F^{2}=0$, $\C_{0}\cdot\F=1$ et $\C_{0}^{2}=-e=\textup{deg}(\G)$.\\
\indent L'application compos\'ee
$${p_{\sur}^{*}\Omega_{\B}^{1}}
\overset{dp_{\sur}}{\longrightarrow}\Omega_{\sur}^{1}
\overset{\eta_{\sur}}{\longrightarrow}\lb_{\sur}^{-1}$$
est nulle sur les fibres g\'en\'erales de $p_{\sur}$
et donc identiquement nulle. Le feuilletage $(\eta_{\sur},\lb_{\sur})$ 
est donc vertical relativement \`a $p_{\sur}$, autrement dit, il existe une factorisation\\
\centerline{
\xymatrix{
{\Omega_{\sur}^{1}} \ar[r]\ar[d]^-{\eta_{\sur}} 
& {\Omega_{\sur/\B}^{1}} \ar@/^/[dl]\\
\lb_{\sur}^{-1}
}
}
\vspace{-1mm}\\
ou encore, des inclusions $\lb_{\sur}\subset\T_{\sur/\B}\subset\T_{\sur}$. Il 
existe donc un diviseur effectif $\D_{\sur}\subset\sur$ 
tel que $\T_{\sur/\B}=\lb_{\sur}(\D_{\sur})$ et $\D_{\sur}\cdot\F=1$.\\
\indent Soit
$\C_{0}+b\F$ la classe de $\D_{\sur}$ dans $\textup{Num}(\sur)$. 
La classe de $\lb_{\sur}$ dans $\textup{Num}(\sur)$ est donc $\C_{0}+(e-b)\F$ 
et
$$\lb_{\sur}\cdot\D_{\sur}=(\C_{0}+(e-b)\F)\cdot(\C_{0}+b\F)=-e+b+(e-b)=0.$$ 
Si $\D_{\sur}^{h}$ est la composante horizontale de $\D_{\sur}$
et $\D_{\sur}^{v}$ sa composante verticale, alors
$$0=\lb_{\sur}\cdot\D_{\sur}=\lb_{\sur}\cdot\D_{\sur}^{h}+\lb_{\sur}\cdot\D_{\sur}^{v}$$ 
et, puisque $\lb_{\sur}$ est num\'eriquement effectif, 
$$\lb_{\sur}\cdot\D_{\sur}^{h}=\lb_{\sur}\cdot\D_{\sur}^{v}=0$$ 
et
$$\D_{\sur}^{v}=0.$$ 
Le lieu singulier $\D_{\sur}$ du feuilletage $(\eta_{\sur},\lb_{\sur})$ 
est donc une section de $p_{\sur}$ contract\'ee par $q_{\sur}$. L'unicit\'e d'une telle
section est bien connue.\\
$\ $
\newline
\noindent $2^{\textup{\`eme}}$ cas : $\lb_{\sur}\cdot\F=2$.$-$Il 
existe une r\'esolution \textit{minimale} et
\textit{\'equivariante} $d : \sur_{1}\longrightarrow\sur$
des singularit\'es de $\sur$ \re{BuWa74}{Proposition 1.2}. 
Notons $\eta_{1} : \Omega_{\sur_{1}}^{1}\longrightarrow d^{*}\lb_{\sur}^{-1}$ 
le rel\`evement de $\eta_{\sur} : \Omega_{\sur}^{1}
\longrightarrow\lb_{\sur}^{-1}$ 
\`a $\sur_{1}$ (\textit{voir} Lemme 1.2).\\
\indent Le genre de
$\B$ est $\ge 1$ et il n'y a donc pas de ($-1$)-courbe sur $\sur_{1}$ dominant $\B$.
Supposons que la surface $\sur_{1}$ ne soit pas minimale.  
Son mod\`ele minimal est une surface g\'eom\'etriquement r\'egl\'ee
$\sur_{k}$ et s'obtient \`a partir de $\sur_{1}$ en contractant des ($-1$)-courbes
dans les fibres de $p_{1}$\\
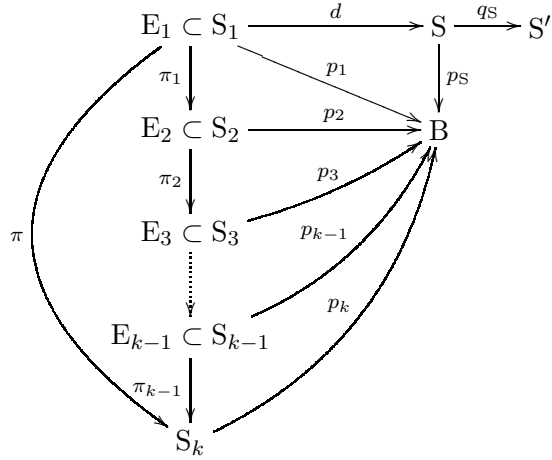
\centerline{
\xymatrix{
{\E_{1}\subset\sur_{1}} \ar[rr]^-{d}\ar[d]_-{\pi_{1}} \ar[drr]^-{p_{1}}
\ar@/_5pc/[dddd]_-{\pi}
 & & {\sur} \ar[d]^-{p_{\sur}}\ar[r]^-{q_{\sur}} & \sur'\\
{\E_{2}\subset\sur_{2}} \ar[rr]^-{p_{2}} \ar[d]_-{\pi_{2}} &  & {\B}\\
{\E_{3}\subset\sur_{3}} \ar@/_/[urr]^-{p_{3}} \ar@{.>}[d]\\
{\E_{k-1}\subset\sur_{k-1}} \ar[d]_-{\pi_{k-1}}\ar@/_1.2pc/[uurr]^-{p_{k-1}}\\
{\sur_{k}} \ar@/_1.7pc/[uuurr]^-{p_{k}}
}
}
\vspace{-1mm}\\
Le morphisme $\pi_{i}$ est l'\'eclatement d'un point de $\sur_{i+1}$, de lieu
exceptionnel $\E_{i}\subset\sur_{i}$. Notons $\lb_{1}=d^{*}\lb_{S}$. 
Le fibr\'e $\lb_{1}$ est num\'eriquement effectif et 
$\lb_{1}\cdot\F=\lb_{\sur}\cdot\F=2$. 
Il existe un fibr\'e $\lb_{2}$ sur $\sur_{2}$ et un entier $j_{1}\in\ZZ$ tels que
$\lb_{1}=\pi_{1}^{*}\lb_{2}(j_{1}\E_{1})$. Comme $\lb_{1}$ est num\'eriquement 
effectif, $j_{1}\le 0$ et $\lb_{2}$ est \'egalement num\'eriquement effectif. 
L'inclusion
$\T_{\sur_{1}}\longrightarrow\pi_{1}^{*}\T_{\sur_{2}}$ donne les inclusions
$$0\longrightarrow\T_{\sur_{1}}\otimes\lb_{1}^{-1}
\longrightarrow\pi_{1}^{*}(\T_{\sur_{2}}\otimes\lb_{2}^{-1})(-j_{1}\E_{1})$$
et
$$0\longrightarrow{\pi_{1}}_{*}(\T_{\sur_{1}}\otimes\lb_{1}^{-1})
\longrightarrow\T_{\sur_{2}}\otimes\lb_{2}^{-1}$$
puisque 
${\pi_{1}}_{*}\mathcal{O}_{\sur_{1}}(-j_{1}\E_{1})=\mathcal{O}_{\sur_{2}}$.
En particulier, $\coh^{0}(\sur_{1},\T_{\sur_{1}}\otimes\lb_{1}^{-1})\subset
\coh^{0}(\sur_{2},\T_{\sur_{2}}\otimes\lb_{2}^{-1})$. L'argument pr\'ec\'edent montre
qu'il existe, pour tout $2\le i\le k$, des fibr\'es num\'eriquement effectifs $\lb_{i}$ 
sur
$\sur_{i}$ et des entiers $j_{i}\le 0$ tels que 
$$\lb_{i}=\pi^{*}\lb_{i+1}(j_{i}\E_{i})$$
et 
$$\coh^{0}(\sur_{i-1},\T_{\sur_{i-1}}\otimes\lb_{i-1}^{-1})\subset
\coh^{0}(\sur_{i},\T_{\sur_{i}}\otimes\lb_{i}^{-1}).$$ 
D'o\`u
$$\lb_{1}=\pi^{*}\lb_{k}(-\E)$$ 
o\`u $\E$ est un diviseur effectif et $\pi$-exceptionnel.
Soit enfin $\eta_{k}$ l'image de 
$\eta_{1}\in\coh^{0}(\sur_{1},\T_{\sur_{1}}\otimes\lb_{1}^{-1})$ 
dans $\coh^{0}(\sur_{k},\T_{\sur_{k}}\otimes\lb_{k}^{-1}).$\\
\indent Soit $\G$ un fibr\'e de rang 2
sur $\B$, \textit{normalis\'e}, de sorte
que $\sur_{k}=\PP_{\B}(\G)$. Notons $\C_{0}\subset\sur_{k}$ une section de $p_{k}$
telle que $\mathcal{O}_{\sur_{k}}(\C_{0})$ soit isomorphe au fibr\'e tautologique
$\mathcal{O}_{\sur_{k}}(1)$. Soit $e=-\textup{deg}(\G)$.\\
\indent L'application compos\'ee
$${p_{k}^{*}\Omega_{\B}^{1}}
\overset{dp_{k}}{\longrightarrow}\Omega_{\sur_{k}}^{1}
\overset{\eta_{k}}{\longrightarrow}\lb_{k}^{-1}$$
est nulle sur les fibres g\'en\'erales de $p_{k}$
et donc identiquement nulle. Le feuilletage $(\eta_{k},\lb_{k})$ 
est donc vertical relativement \`a $p_{k}$, autrement dit, il existe une factorisation\\
\centerline{
\xymatrix{
{\Omega_{\sur_{k}}^{1}} \ar[r]\ar[d]^-{\eta_{k}} 
& {\Omega_{\sur_{k}/\B}^{1}} \ar@/^/[dl]\\
\lb_{k}^{-1}
}
}
\vspace{-1mm}\\
ou encore, des inclusions $\lb_{k}\subset\T_{\sur_{k}/\B}\subset\T_{\sur_{k}}$. Il 
existe donc un diviseur effectif $\D_{k}\subset\sur_{k}$ 
tel que $\T_{\sur_{k}/\B}=\lb_{k}(\D_{k})$. Enfin, $\lb_{k}\cdot\F=2$.\\
\indent Le fibr\'e tangent relatif $\T_{\sur_{k}/\B}$ est de
degr\'e 2 sur la fibre $\F$ et $\D_{k}$ est donc vertical relativement \`a $p_{k}$.
La classe de $\T_{\sur_{k}/\B}$ dans $\textup{Num}(\sur_{k})$ est $2\C_{0}+e\F$ et la classe
de $\lb_{k}$ est donc $2\C_{0}+(e-\alpha)\F$ o\`u $\alpha=\D_{k}\cdot\C_{0}\ge 0$. 
Le nombre $\lb_{k}^{2}=-4\alpha$ est donc $\le 0$ d'une part et $\ge 0$ d'autre part
puisque $\lb_{k}$ est num\'eriquement effectif, et donc nul. D'o\`u
$$\D_{k}=0\text{ et }
\lb_{k}\equiv 2\C_{0}+e\F.$$ 
\indent La relation $\lb_{1}=\pi^{*}\lb_{k}(-\E)$ donne 
$\lb_{1}^{2}=\lb_{k}^{2}+\E^{2}=\E^{2}$. Le premier de ces deux nombres est $\ge 0$ 
puisque $\lb_{1}$ est num\'eriquement effectif et le second est $\le 0$ puisque
$\E$ est $\pi$-exceptionnel. D'o\`u 
$$\E=0\text{ et }\lb_{1}\cdot\E_{1}=\pi^{*}\lb_{k}\cdot\E_{1}=0.$$ La 
$(-1)$-courbe $\E_{1}$ est donc contract\'ee par $q_{\sur}\circ d$ mais ne l'est pas 
par $d$, par minimalit\'e de la r\'esolution $d : \sur_{1}\longrightarrow\sur^{n}$.
La courbe $d(\E_{1})$ n'est pas contract\'ee par $q_{\sur}$ 
car elle est verticale relativement \`a $p_{1}$. La 
surface $\sur_{1}$ est donc minimale et, par ce qui pr\'ec\'ede, $\lb_{1}$ est isomorphe
au fibr\'e tangent relatif $\T_{\sur_{1}/\B}$.\\ 
\indent Supposons qu'il existe une courbe irr\'eductible
$\C\subset\sur_{1}$ qui soit contract\'ee par $q_{\sur}\circ d$. Soit 
$a\C_{0}+b\F$ sa classe dans $\textup{Num}(\sur_{1})$,
$$\lb_{1}\cdot\C=2b-ae\text{ et }\C^{2}=a(2b-ae).$$ 
Le premier de ces deux nombres doit 
\^etre nul si $\C$ est contract\'ee par $q_{\sur}\circ d$ et le second strictement
n\'egatif, ce qui est absurde. La surface $\sur$ est donc lisse, 
et le morphisme $q_{\sur}$ est fini. Enfin, le feuilletage 
$(\eta_{\sur},\lb_{\sur})$ sur $\sur$ est r\'egulier.\qed
\sect{Fibr\'es projectifs}
\marque Soit $\Y$ une vari\'et\'e alg\'ebrique d\'efinie
sur un corps alg\'ebriquement clos. 
Soit $\G$ un fibr\'e vectoriel de rang 2 sur $\Y$, extension d'un fibr\'e en droites
$\M$ par le fibr\'e trivial $\mathcal{O}_{\Y}$,
$$0\longrightarrow\mathcal{O}_{\Y}\longrightarrow\G\longrightarrow\M
\longrightarrow 0.$$
Soit $\alpha$ la classe dans $\coh^{1}(\Sigma,\M^{-1})$ de cette extension.
Soient $\Z:=\PP_{\Y}(\G)$ et $p:\Z\longrightarrow\Y$ le morphisme naturel.
Soit $\Sigma$ la section de $p$ correspondant 
au quotient inversible $\M$ de $\G$. L'id\'eal 
de $\Sigma$ dans $\Z$ s'identifie au fibr\'e $\mathcal{O}_{\Z}(-1)$.\\
\indent Soit $2\Sigma$ le deuxi\`eme voisinage
infinit\'esimal de $\Sigma$ dans $\Z$, c'est-\`a-dire, le sous-sch\'ema ferm\'e
de $\Z$ d\'efini par l'id\'eal $\I_{\Sigma/\Z}^{2}$. 
\lemme{Il existe un isomorphisme 
$\textup{Pic}(2\Sigma)\simeq\textup{Pic}(\Sigma)\oplus\coh^{1}(\Sigma,\M^{-1}).$}
\textit{D\'emonstration}.$-$Le morphisme $p$ induit
une r\'etraction de l'inclusion $\Sigma\subset 2\Sigma$. La suite exacte courte
$$0\longrightarrow\M^{-1}\simeq\I_{\Sigma/\Z}/\I_{\Sigma/\Z}^{2}
\longrightarrow\mathcal{O}^{\times}_{2\Sigma}
\longrightarrow\mathcal{O}^{\times}_{\Sigma}\longrightarrow 0$$
donne la suite exacte de cohomologie
$$0\longrightarrow\coh^{1}(\Sigma,\M^{-1})\longrightarrow
\coh^{1}(\Sigma,\mathcal{O}^{\times}_{2\Sigma})\longrightarrow
\coh^{1}(\Sigma,\mathcal{O}^{\times}_{\Sigma}).$$
La r\'etraction $2\Sigma\longrightarrow\Sigma$ induit
en cohomologie une section de $\coh^{1}(\Sigma,\mathcal{O}^{\times}_{2\Sigma})\longrightarrow
\coh^{1}(\Sigma,\mathcal{O}^{\times}_{\Sigma})$ et, finalement,
$$\textup{Pic}(2\Sigma)\simeq\textup{Pic}(\Sigma)\oplus\coh^{1}(\Sigma,\M^{-1}).$$\qed\\
\indent Le point suivant a \'et\'e observ\'e par Schr\"{o}er.
\lemmeref{Sc00}{Example 5.6}{La classe du fibr\'e 
$\mathcal{O}_{\Z}(1)\otimes p^{*}\M^{-1}$ 
dans $\textup{Pic}(2\Sigma)$ est $(0,\alpha)$.}
\textit{D\'emonstration}.$-$La classe du fibr\'e $\mathcal{O}_{\Z}(1)\otimes p^{*}\M^{-1}$
dans $\textup{Pic}(\Sigma)$ est nulle. Soit $\U_{i}$ un recouvrement par des ouverts affines
de $\Sigma\simeq\Y$ trivialisant le fibr\'e $\M$. Soit $\M_{i}$ la restriction de $\M$
\`a $\U_{i}$ et $s_{i}\in\coh^{0}(\U_{i},\M_{i})$ une section partout non nulle.
Soient $\psi_{i}:\M_{i}\longrightarrow\mathcal{O}_{i}$ la trivialisation de $\M_{i}$
induite par $s_{i}$ et $f_{ij}=\psi_{i}\circ\psi_{j}^{-1}$ les fonctions
de transitions du fibr\'e $\M_{i}$ sur les ouverts
$\U_{ij}:=\U_{i}\cap\U_{j}$. Soit $\G_{i}$ la restriction de $\G$ \`a 
$\U_{i}$ et fixons un isomorphisme 
$\varphi_{i}:\G_{i}\longrightarrow\mathcal{O}_{i}\oplus\M_{i}$  
rendant le diagramme suivant commutatif\\
\centerline{
\xymatrix{
0 \ar[r] & {\mathcal{O}_{i}} \ar[r] \ar@{=}[d]
& {\G_{i}} \ar[r] \ar[d]^-{\varphi_{i}}
& {\M_{i}} \ar[r] \ar[d]^-{\psi_{i}}
& 0\\
0 \ar[r] & {\mathcal{O}_{i}} \ar[r]
& {\mathcal{O}_{i}\oplus\mathcal{O}_{i}} \ar[r]
& {\mathcal{O}_{i}} \ar[r]
& 0 
}
}
\vspace{-1mm}\\
Les fonctions de transitions $\varphi_{i}\circ\varphi_{j}^{-1}$, sur les ouverts
$\U_{ij}$, du fibr\'e $\G$ sont de la forme 
$$
\left(
\begin{array}{cc}
1 & g_{ij}\\
0 & f_{ij}
\end{array}
\right)
$$
\vspace{-1mm}\\
Soit $s_{i}^{*}$ la section duale
de $s_{i}$ sur $\U_{i}$ et posons, pour tout $(i,j)$, 
$\alpha_{ij}=g_{ij}{s_{j}^{*}}_{|U_{ij}}=\frac{g_{ij}}{f_{ij}}{s_{i}^{*}}_{|U_{ij}}$. 
La classe $\alpha\in\coh^{1}(\Sigma,\M^{-1})$
est alors represent\'ee par le cocycle $(\alpha_{ij})$.\\
\indent Il reste \`a faire le calcul des fonctions de transitions du fibr\'e
$\mathcal{O}_{\Z}(1)$ au voisinage de la section $\Sigma$. Posons 
$\Z_{i}=p^{-1}(\U_{i})$ et $\Z_{ij}:=\Z_{i}\cap\Z_{j}$.
Notons $\hat{s}_{i}=\varphi_{i}^{-1}(0,1)$ la section
de $\mathcal{O}_{\Z}(1)$ au dessus de l'ouvert $\Z_{i}$. 
Sur $\Z_{i}\cap\{\hat{s}_{i}\neq 0\}$, $\varphi_{i}^{-1}(1,0)=t_{i}\hat{s}_{i}$,
o\`u $t_{i}$ est une fonction r\'eguli\`ere sur l'ouvert consid\'er\'e.
Sur $\Z_{ij}\cap\{\hat{s}_{i}\neq 0\}$,
$$\hat{s}_{j}=(p^{*}f_{ij}+p^{*}g_{ij}t_{i})\hat{s}_{i}.$$
Notons que les ouverts 
$\Z_{i}\cap\{\hat{s}_{i}\neq 0\}$ de $\Z$ recouvrent $\Sigma$.
Les fonctions de transitions $h_{ij}$ du fibr\'e tautologique
sur les ouverts $(\Z_{i}\cap\{\hat{s}_{i}\neq 0\})\cap(\Z_{j}\cap\{\hat{s}_{j}\neq 0\})$
sont donc donn\'ees par la formule
$$\frac{h_{ij}}{p^{*}f_{ij}}=1+
\frac{p^{*}g_{ij}}{p^{*}f_{ij}}t_{i}.$$
Il reste \`a remarquer que 
le ferm\'e $\Sigma$ est d\'efini, sur l'ouvert $\Z_{i}\cap\{\hat{s}_{i}\neq 0\}$,
par l'annulation de la fonction $t_{i}$.\qed
\sect{D\'emonstration des r\'esultats}
\marque Soient $\X$ une vari\'et\'e projective \textit{normale}, 
de dimension $n\ge 2$, d\'efinie sur un corps $\corps$ alg\'ebriquement clos 
de caract\'eristique nulle, et $(\eta,\lb)$ un feuilletage en courbes sur $\X$.
Supposons $\lb\cdot\C>0$ pour toute courbe $\C\subset\X$.\\
$\ $
\newline
\marque\textit{Existence de feuilles alg\'ebriques}.$-$Notons $\X_{s}$ le lieu 
singulier de $\X$ et $\X_{ns}$ le compl\'ementaire de $\X_{s}$ dans $\X$. Soit
$\D_{sat}^{ns}$ le diviseur des composantes de codimension 1, 
compt\'ees avec multiplicit\'es, 
de $\Z(\eta,\lb)\cap\X_{ns}$ et $\D_{sat}$ l'adh\'erence de son \textit{support} dans $\X$. 
Notons 
$\eta_{sat}^{ns} : \Omega_{\X_{ns}}^{1}\longrightarrow{\lb^{-1}}_{|\X_{ns}}(- \D_{sat}^{ns})$ 
le feuilletage induit : le lieu 
$\Z(\eta_{sat}^{ns},{\lb}_{\X_{ns}}(\D_{sat}^{ns}))$ des
z\'eros de $\eta_{sat}^{ns}$ est de codimension au moins 2 dans $\X_{ns}$. Notons enfin
$\X_{r}$ le compl\'ementaire dans $\X$ de 
$\overline{{\Z(\eta_{sat}^{ns},{\lb}_{|\X_{ns}}(\D_{sat}^{ns}))}_{red}}\cup 
\overline{{\Z(\eta,\lb)}_{red}-\D_{sat}}$. Remarquons enfin que le ferm\'e $\X-\X_{r}$
est de codimension au moins 2 dans $\X$.\\
\indent Le ferm\'e $\X-(\X_{ns}\cap\X_{r})$ est de codimension au moins 2 dans $\X$.
Il existe donc une courbe projective lisse
$\B\subset\X_{ns}\cap\X_{r}$ telle que
\begin{enumerate}
\item[$\bullet$] $\B$ n'est pas une feuille du feuilletage
$(\eta_{sat}^{ns},{\lb}_{|\X_{ns}}(\D_{sat}^{ns}))$,
\item[$\bullet$] $\B\not\subset\D_{sat}$.
\end{enumerate} 
Notons $\Gamma_{\B}\subset\B\times(\X_{ns}\cap\X_{r})\subset\B\times\X$ 
le graphe de $\B\subset\X$. Notons $p$ et $q$ les projections de $\B\times\X$ sur 
$\B$ et $\X$ respectivement, $p_{ns}$ et $q_{ns}$ leurs restrictions \`a 
$\B\times\X_{ns}$. 
Le feuilletage $(\eta_{sat}^{ns},{\lb}_{|\X_{ns}}(\D_{sat}^{ns}))$ induit un 
feuilletage sur $\B\times\X_{ns}$, 
$$\Omega_{\B\times\X_{ns}}^{1}
=p_{ns}^{*}\Omega_{\B}^{1}\oplus q_{ns}^{*}\Omega_{\X_{ns}}^{1}
\longrightarrow q_{ns}^{*}\Omega_{\X_{ns}}^{1}\longrightarrow 
q_{ns}^{*} {\lb^{-1}}_{|\X_{ns}}(- \D_{sat}^{ns}).$$ 
Ce feuilletage est partout 
non nul le long de $\Gamma_{\B}$ et $\Gamma_{\B}$ ne lui est tangente en aucun 
point. Par le th\'eor\`eme de Frobenius formel \res{Mi87}, il existe un sous-sch\'ema formel 
$\widehat{\V}$ lisse de dimension 2 du compl\'et\'e formel $\widehat{\B\times\X}_{ns}$ de 
$\B\times\X_{ns}$ le long de $\Gamma_{\B}$ admettant $\Gamma_{\B}$ comme
sch\'ema de d\'efinition. Le fibr\'e normal $\N_{\Gamma_{B}/\widehat{\V}}$ de 
$\Gamma_{\B}$ dans $\widehat{\V}$ est naturellement isomorphe \`a la restriction
de $q_{ns}^{*}{\lb}_{|\X_{ns}}(\D_{sat}^{ns})$ 
\`a $\Gamma_{\B}$ et en particulier ample. Le sch\'ema formel $\widehat{\V}$ est donc 
alg\'ebrique \re{Bo01}{Theorem 3.5}, autrement dit, il existe une surface $\sur\subset\B\times\X$, 
irr\'eductible
et r\'eduite, telle que $\Gamma_{\B}\subset\sur$ et $\widehat{\V}\subset\widehat{\sur}$ o\`u 
$\widehat{\sur}$ est le compl\'et\'e formel de $\sur$ le long de $\Gamma_{\B}$.\\
\indent Notons $q_{\sur}$ la restriction de $q$ \`a $\sur$. Les 
germes de feuilles du feuilletage $(\eta_{\sur},q_{\sur}^{*}\lb)$, 
aux points de $\Gamma_{\B}\subset\sur$, sont donc alg\'ebriques et 
les images desdites feuilles par $q_{\sur}$
sont des feuilles du feuilletage
$(\eta_{sat}^{ns},{\lb}_{|\X_{ns}}(\D_{sat}^{ns}))$ sur $\X_{ns}$ et, 
sauf \'eventuellement un nombre fini, des feuilles du feuilletage
$(\eta,\lb)$ sur $\X$, puisque $\B\not\subset\D_{sat}$.
\rem Supposons $\coh^{0}(\X,\lb)\neq 0$. Soit $\Y\in|\lb|$ et soit $\G$
la composante neutre du groupe des automorphismes de $\X$ fixant
$\Y$ points par points. Le groupe $\G$ est alg\'ebrique et son alg\`ebre de Lie
$\coh^{0}(\X,\der{\X}(-\Y))$ est, par hypoth\`ese, de 
dimension $>0$. Soit $\G_{0}$ un sous-groupe ferm\'e, 
affine et connexe de $\G$, distingu\'e dans $\G$ et maximal pour ces propri\'et\'es.
Le quotient $\G/\G_{0}$ est alors une vari\'et\'e ab\'elienne par 
le th\'eor\`eme de Chevalley. Supposons le groupe $\G_{0}$ r\'eduit 
\`a l'identit\'e. Le diviseur $\Y$ rencontre, par hypoth\`ese, toutes les courbes 
trac\'ees sur $\X$. Les orbites sous $\G$ non r\'eduites \`a un point sont 
compactes et rencontrent donc le lieu des points fixes sous $\G$, ce qui est absurde.
Il existe donc un sous-groupe affine de $\G$ isomorphe 
\`a $\GG_{a}$ ou $\GG_{m}$ et les orbites sous ce sous-groupe sont les feuilles
d'un feuilletage en courbes, \'eventuellement distinct du feuilletage consid\'er\'e.\\
$\ $
\newline
\marque\textit{La vari\'et\'e des feuilles}.$-$Soit $\bar{\Y}$ l'adh\'erence du lieu des points 
du sch\'ema de Hilbert de $\X$ correspondants aux feuilles du feuilletage 
$(\eta,\lb)$ et soit $\bar{\Z}\subset\bar{\Y}\times\X$ la famille universelle. Notons $p$ et 
$q$ les projections de $\bar{\Y}\times\X$ sur $\bar{\Y}$ et $\X$ respectivement, et 
$p_{\bar{\Z}}$, $q_{\bar{\Z}}$ 
leurs restrictions \`a $\bar{\Z}$. Le sch\'ema $\bar{\Y}$ a 
un nombre au plus d\'enombrable de composantes irr\'eductibles, toutes de type fini
sur le corps de base infini $\corps$. Le lieu des feuilles du 
feuilletage $(\eta,\lb)$ rencontrant 
l'ouvert $\X_{ns}\cap\X_{r}-\D_{sat}$ n'\'etant pas r\'eunion au plus d\'enombrable
de ferm\'es propres de $\X$, il existe une composante irr\'eductible $\Y$
de $\bar{\Y}$ telle que la restriction de $q_{\bar{\Z}}$ \`a la famille universelle 
correspondante soit un morphisme dominant. Soient $\Z\subset\Y\times\X$ la famille universelle, 
$p_{\Z}$ et $q_{\Z}$, les morphismes naturels vers $\Y$ et $\X$ respectivement\\
\centerline{
\xymatrix{
{\Z} \ar[r]^-{q_{\Z}}\ar[d]_-{p_{\Z}} & {\X}\\
{\Y}
}
}
\vspace{-1mm}\\
\lemme{Soient $\Z$ et $\Y$ deux sch\'emas de type fini sur un corps $\corps$ et
$p:\Z\longrightarrow\Y$ un morphisme sur $\corps$, propre et plat. Si $\Y$ est int\`egre et
si l'une des fibres de $p$ est g\'eom\'etriquement int\`egre, alors $\Z$ est int\`egre.}
\noindent\textit{D\'emonstration}.$-$L'irr\'eductibilit\'e de $\Z$ est bien connue.
Il reste \`a v\'erifier que $\Z$ est r\'eduit. Le lieu des points $y\in\Y$ tel 
que $p^{-1}(y)$ soit g\'eom\'etriquement int\`egre est ouvert dans $\Y$ et non vide par
hypoth\`ese. La fibre g\'en\'erique est donc g\'eom\'etriquement int\`egre. 
Soient $\textup{Spec}(\B)\subset\Z$ et $\textup{Spec}(\A)\subset\Z$
deux ouverts affines non vides tels que $p(\textup{Spec}(\B))\subset\textup{Spec}(\A)$.
L'anneau $\B$ est plat sur $\A$, par hypoth\`ese. Si $b\in\B$ est nilpotent alors,
comme la fibre g\'en\'erique est int\`egre, il existe $a\in\A$ non nul, tel 
que $ab=0$. L'\'el\'ement $a$ est non diviseur de z\'ero dans $\A$ et donc dans 
$\B$, et $b=0$.\qed\\
\newline
\indent La famille universelle $\Z$ est donc int\`egre par le lemme
pr\'ec\'edent. Le feuilletage
$(\eta,\lb)$ sur $\X$ induit un 
feuilletage sur $\Y\times\X$ 
$$\Omega_{\Y\times\X}^{1}
=p^{*}\Omega_{\Y}^{1}\oplus q^{*}\Omega_{\X}^{1}
\longrightarrow q^{*}\Omega_{\X}^{1}\longrightarrow 
q^{*}\lb^{-1}.$$ 
\indent L'application compos\'ee
$$\I_{\Z/\Y\times\X}/ \I_{\Z/\Y\times\X}^{2}\longrightarrow
{\Omega_{\Y\times\X}^{1}}_{|\Z}\longrightarrow q^{*}{\lb^{-1}}_{|\Z}$$
est nulle
sur une partie dense de $\Z$ et donc identiquement nulle puisque le faisceau
$q^{*}{\lb^{-1}}_{|\Z}$ est sans torsion. Il existe donc une factorisation
$\eta_{\Z} : \Omega_{\Z}^{1}\longrightarrow q^{*}\lb^{-1}_{|\Z}$ 
et un diagramme commutatif\\
\centerline{
\xymatrix{
{{\Omega_{\Y\times\X}^{1}}_{|\Z}} \ar[r]\ar[r]\ar[d] 
& {\Omega_{\Z}^{1}}\ar@/^/[dl]^-{\eta_{\Z}}   \\  
{q^{*}{\lb^{-1}}_{|\Z}}
}
}
\vspace{-1mm}\\
Le lieu singulier $\Z(\eta_{\Z},q_{\Z}^{*}\lb)$ est intersection
sch\'ematique des ferm\'es $\Y\times\Z(\eta,\lb)$ et $\Z$ de $\Y\times\X$. 
L'application compos\'ee
$$p_{\Z}^{*}\Omega_{\Y}^{1}\longrightarrow\Omega_{\Z}^{1}
\longrightarrow q_{\Z}^{*}\lb^{-1}$$
est nulle
sur les fibres g\'en\'erales de $p_{\Z}$ et donc identiquement nulle. 
Le feuilletage $(\eta_{\Z},q_{\Z}^{*}\lb)$ 
est donc vertical relativement \`a $p_{\Z}$, autrement dit, il existe une factorisation\\
\centerline{
\xymatrix{
{\Omega_{\Z}^{1}} \ar[r]\ar[d]^-{\eta_{\Z}} 
& {\Omega_{\Z/\Y}^{1}} \ar@/^/[dl]\\
q_{\Z}^{*}\lb^{-1}
}
}
\vspace{-1mm}\\
Le morphisme $q_{\Z}$ est donc birationnel.\\
$\ $
\newline
\marque\textit{R\'eduction au cas des surfaces}.$-$Soient $\B$ une courbe projective, 
lisse et connexe, de genre $g(\B)>0$ et $\B\longrightarrow\Y$ un morphisme non constant, 
dont l'image
rencontre le lieu des feuilles, de sorte que $\eta_{\Z}$ est non identiquement 
nulle le long d'au moins une fibre de $p_{\Z}$ au dessus de $\B$. 
Soit $\sur:=\B\times_{\Y}\Z\subset\B\times\X$ ; $\sur$ est
une surface irr\'eductible et r\'eduite (\textit{voir} Lemme 4.5). 
Soient $p_{S}$ et $q_{S}$ les restrictions \`a $\sur$ des projections de 
$\B\times\X$ sur $\B$ et $\X$ respectivement. Le feuilletage
$(\eta,\lb)$ induit, comme ci-dessus, un feuilletage en courbes
$\eta_{\sur} : \Omega_{\sur}^{1}\longrightarrow q_{\sur}^{*}\lb^{-1}$.
Le lieu singulier $\Z(\eta_{\sur},q_{\sur}^{*}\lb)$ est intersection
sch\'ematique des ferm\'es $\B\times\Z(\eta,\lb)$ et $\sur$ de $\B\times\X$.
Soit $\lb_{\sur}$ le fibr\'e $q_{\sur}^{*}\lb$. Soit 
$n : \sur^{n}\longrightarrow\sur$ la normalisation de $\sur$ et
$\eta_{n} : \Omega_{\sur^{n}}^{1}\longrightarrow n^{*}\lb_{\sur}^{-1}$
le rel\`evement de 
$\eta_{\sur} : \Omega_{\sur}^{1}\longrightarrow \lb_{\sur}^{-1}$
\`a $\sur^{n}$ (\textit{voir} Lemme 1.2). Soit $\sur':=q_{\sur}(\sur)\subset\X$ et soit
$q_{\sur^{n}}:=q_{\sur}\circ n$.
Le fibr\'e $\lb_{\sur^{n}}$ num\'eriquement effectif et
une courbe $\C\subset\sur^{n}$ est contract\'ee par $q_{\sur^{n}}$ si et seulement si
$\lb_{\sur^{n}}\cdot\C=0$. Enfin, les courbes trac\'ees sur $\sur^{n}$ et
contract\'ees par $q_{\sur^{n}}$ sont horizontales relativement \`a 
$p_{\sur^{n}}:=p_{\sur}\circ n$. Les hypoth\`eses de la proposition 2.1 sont satisfaites.
La surface $\sur^{n}\longrightarrow\B$ est donc g\'eom\'etriquement
r\'egl\'ee et $n^{*}\lb_{\sur}\cdot\F\in\{1,2\}$, o\`u $\F$ est une 
fibre g\'en\'erale de $p_{\sur^{n}}$.\\
$\ $
\newline
\noindent $1^{er}$ cas : $n^{*}\lb_{\sur}\cdot\F=2$.$-$Le 
morphisme $q_{\sur}\circ n$ est fini sur son image et le feuilletage 
$(\eta_{n},n^{*}\lb_{\sur})$ sur $\sur^{n}$ est r\'egulier (\textit{voir} Proposition 2.1).
Supposons $\Z(\eta,\lb)\neq\emptyset$ et supposons que 
$q_{\sur}(\sur)\cap\Z(\eta,\lb)\neq\emptyset$. 
Le lieu singulier 
$\Z(\eta_{\sur^{n}},n^{*}\lb)$ alors non vide, puisque 
$$n^{-1}((\B\times\Z(\eta,\lb)_{\textup{red}})\cap\sur)
=n^{-1}({\Z(\eta_{\sur},\lb_{\sur})}_{\textup{red}})
\subset {\Z(\eta_{\sur^{n}},n^{*}\lb_{\sur})}_{\textup{red}},$$
ce qui est absurde. Le feuilletage $(\eta,\lb)$ est donc r\'egulier.\\
\newline
\indent Si $\X$ est une surface normale, le morphisme $q_{\Z}$ est birationnel et 
fini : c'est un isomorphisme.\\
$\ $
\newline
\noindent $2^{\text{\`eme}}$ cas : $n^{*}\lb_{\sur}\cdot\F=1$.$-$Les 
fibres du morphisme $p_{\Z}$ sont donc des courbes rationnelles 
irr\'eductibles et g\'en\'eriquement r\'eduites. Les fibres g\'en\'erales sont r\'eduites et
irr\'eductibles.\\
\indent Le lieu singulier $\Z(\eta_{\sur^{n}},n^{*}\lb)$ est irr\'eductible de 
dimension 1 (\textit{voir} Proposition 2.1) et contract\'e par $q_{\Z}$. 
Les feuilles du feuilletage $(\eta,\lb)$ passent donc par un point $x\in\X$.
Supposons $\Z(\eta,\lb)_{\textup{red}}\neq\{x\}$ et soit $z\in\Z(\eta,\lb)_{\textup{red}}$
distinct de $x$. Supposons que $z\in q_{\sur}(\sur)$. La proposition 2.1 nous permet
de supposer que l'une des feuilles du feuilletage
$(\eta,\lb)$ param\'etr\'ees par $\B$ ne passe pas par $z$. Soit $b\in\B$ tel
que $(b,z)\in\sur$. Le point $(b,z)$ est dans le lieu singulier
${\Z(\eta_{\sur},\lb_{\sur})}_{\textup{red}}=(\B\times\Z(\eta,\lb)_{\textup{red}})\cap\sur$
et, comme
$$n^{-1}({\Z(\eta_{\sur},\lb_{\sur})}_{\textup{red}})
\subset {\Z(\eta_{\sur^{n}},n^{*}\lb_{\sur})}_{\textup{red}},$$
$n^{-1}(b,z)\subset{\Z(\eta_{\sur^{n}},n^{*}\lb_{\sur})}_{\textup{red}}$.
Or, ${\Z(\eta_{\sur^{n}},n^{*}\lb_{\sur})}_{\textup{red}}$ est contract\'e par $q_{\sur}$ 
sur le point $x$, ce qui est absurde. La remarque pr\'ec\'edente ne sera pas utilis\'ee
par la suite.\\
\indent Soit $\textup{Hom}_{\textup{bir}}(\PP^1,\X,0\mapsto x)$ le sch\'ema des morphismes
birationnels de $\PP^1$ vers $\X$ appliquant $0\in\PP^{1}$ sur $x\in\X$ et soit
$\textup{Hom}_{\textup{bir}}^{n}(\PP^1,\X,0\mapsto x)$
sa normalisation. Le sous-groupe lin\'eaire $\G\subset\textup{PGL}(2)$ 
des automorphismes de $\PP^{1}$ fixant $0$
agit sur $\textup{Hom}_{\textup{bir}}^{n}(\PP^1,\X)$
et $\textup{Hom}_{\textup{bir}}^{n}(\PP^1,\X,0\mapsto x)\times\PP^1$.
Les quotients g\'eom\'etriques au sens de Mumford existent
et seront respectivement not\'es
$\textup{RatCurves}^{n}(\X,x)$ et 
$\textup{Univ}^{rc}(\X,x)$ \res{Ko96}.\\
\indent Les feuilles du feuilletage $(\eta,\lb)$ rencontrant 
l'ouvert $\X_{ns}\cap\X_{r}$ sont alg\'ebriquement \'equivalentes,
par 4.6 par exemple, et l'une des composantes de 
$\textup{RatCurves}^{n}(\X,x)$, not\'ee $\V$, param\`etre les fibres 
de $p_{\Z}$ au dessus d'un ouvert non vide de $\Y$.\\
\indent Le lemme de cassage de Mori entra\^{\i}ne
d'une part que $\V$ est projective et
d'autre part que le morphisme d'\'evaluation
$\U:=\textup{Univ}^{rc}(\X,x)\times_{\textup{RatCurves}^{n}(\X,x)}\V
\longrightarrow\X$ est 
g\'en\'eriquement fini. Le point g\'en\'erique de $\V$ param\`etre donc
une feuille du feuilletage $(\eta,\lb)$ et, par suite, le morphisme
pr\'ec\'edent est en fait birationnel.
Soient $p_{\U}$ et $q_{\U}$ les morphismes naturels\\ 
\centerline{
\xymatrix{
{\U} \ar[r]^-{q_{\U}}\ar[d]_-{p_{\U}} & {\X}\\
{\V}
}
}
\vspace{-1mm}\\
Le morphisme $p_{\U}$ est lisse \re{Ko96}{Corollary II 2.12}.
Le point, au coeur de la preuve du lemme de cassage, est
le r\'esultat bien connu suivant. Si 
$p : \sur\longrightarrow\B$ est une surface g\'eom\'etriquement r\'egl\'ee et 
$q : \sur\longrightarrow\sur'$
est un morphisme g\'en\'eriquement fini, il existe au plus
une courbe, irr\'eductible et r\'eduite, trac\'ee sur $\sur$, horizontale relativement \`a $p$ et 
contract\'ee par $q$. Les fibres du morphisme $q_{\U}$ sont connexes
par le th\'eor\`eme de Zariski et, en particulier, les composantes irr\'eductibles 
des fibres de $q_{\U}$ sont de dimension $>0$. L'une
des composantes de $\Sigma:=q_{\U}^{-1}(x)$ est de codimension 1 dans $\U$. 
On d\'eduit de la remarque pr\'ec\'edente,
que le ferm\'e $\Sigma$ est irr\'eductible de codimension $1$ 
dans $\U$, et que le morphisme $q_{\U}$ induit un isomorphisme 
de $\U-q_{\U}^{-1}(x)$ sur $\X-\{x\}$.\\
\indent Soit $\B\subset\V$ est une courbe lisse g\'en\'erique.
La surface $p_{\U}^{-1}(\B)$ est la normalisation de son image $\sur$ dans $\B\times\X$ et 
il existe un morphisme $\B\longrightarrow\Y$ tel que 
$\sur$ s'identifie au produit fibr\'e $\Z\times_{\Y}\B$.
La trace de $\Sigma$ sur la surface $p_{\U}^{-1}(\B)$ est donc l'unique section de 
$p_{\U}^{-1}(\B)\longrightarrow\B$ contract\'ee par $q_{\U}$ 
(\textit{voir} Proposition 2.1). Le morphisme $\Sigma\longrightarrow\V$ 
induit par $p_{\U}$ est donc birationnel et fini, c'est un isomorphisme.\\
\indent Soit $\G:={p_{\U}}_{*}\mathcal{O}_{\U}(\Sigma)$. Le sch\'ema 
$\U\longrightarrow\V$ 
s'identifie au fibr\'e projectif $\PP_{\V}(\G)\longrightarrow\V$ et le 
fibr\'e $\mathcal{O}_{\U}(\Sigma)$ au fibr\'e tautologique 
$\mathcal{O}_{\U}(1)$. La section $\Sigma$ de $p_{\U}$ correspond 
\`a un quotient inversible $\M$ de $\G$ sur $\U$. Le fibr\'e $\F$ 
est une extension de $\M$ par le fibr\'e trivial $\mathcal{O}_{\V}$
$$0\longrightarrow\mathcal{O}_{\V}\longrightarrow\G\longrightarrow\M
\longrightarrow 0.$$
Soit $\alpha$ la classe dans $\coh^{1}(\Sigma,\M^{-1})$ de cette extension.
Le fibr\'e $q_{\U}^{*}\lb$ est de degr\'e 1 sur les fibres de $p_{U}$
et sa restriction \`a $\Sigma$ est triviale. Il est donc isomorphe au
fibr\'e $\mathcal{O}_{\U}(1)\otimes p_{\U}^{*}\M^{-1}$.\\
\indent Soit $2\Sigma$ le deuxi\`eme voisinage
infinit\'esimal de $\Sigma$ dans $\U$.
La classe du fibr\'e $\mathcal{O}_{\U}(1)\otimes p_{\U}^{*}\M^{-1}$ 
dans $\textup{Pic}(2\Sigma)\simeq\textup{Pic}(\Sigma)\oplus\coh^{1}(\Sigma,\M^{-1})$ 
est $(0,\alpha)$ (\textit{voir} Lemme 3.3) d'une part et nulle d'autre part. 
L'extension d\'efinie par $\alpha$ est donc scind\'ee et le morphisme 
$q_{U}$ s'identifie \`a la contraction de la
section de $\PP_{\V}(\mathcal{O}_{\V}\oplus\M)\longrightarrow\V$ de fibr\'e
normal $\M^{-1}$.
La section de fibr\'e normal trivial ne rencontre pas la section
pr\'ec\'edente, et le diviseur associ\'e dans $\X$, not\'e 
$\mathcal{O}_{\X}(\V)$, s'identifie au fibr\'e $\lb$.\qed

$\ $
\newline
\noindent St\'ephane \textsc{Druel}, \textsc{Institut Fourier}, UMR 5582
du CNRS,
Universit\'e Joseph Fourier, BP 74, 38402 Saint Martin d'H\`eres,
France.\\
e-mail : \texttt{druel@ujf-grenoble.fr}

\end{document}